\pdfoutput=1
\RequirePackage{ifpdf}
\ifpdf 
\documentclass[pdftex]{sigma}
\else
\documentclass{sigma}
\fi

\begin{document}

\allowdisplaybreaks

\renewcommand{\PaperNumber}{007}

\FirstPageHeading

\ShortArticleName{On the Dimension of the Group of Projective Transformations}

\ArticleName{On the Dimension of the Group\\ of Projective Transformations\\ of Closed Randers  and Riemannian Manifolds}

\Author{Vladimir S. MATVEEV}

\AuthorNameForHeading{V.S.~Matveev}

\Address{Institute of Mathematics,  Friedrich-Schiller-Universit\"at Jena,  07737 Jena, Germany}
\Email{\href{mailto:vladimir.matveev@uni-jena.de}{vladimir.matveev@uni-jena.de}}
\URLaddress{\url{http://users.minet.uni-jena.de/~matveev/}}

\ArticleDates{Received January 18, 2012, in f\/inal form February 21, 2012; Published online February 23, 2012}

\Abstract{We construct a counterexample to Theorem~2 of [Raf\/ie-Rad~M., Rezaei~B., {\it SIGMA} {\bf 7} (2011), 085, 12~pages].}


\Classification{ 53C60; 53B40; 53A20}

\section{Introduction}

A Randers  metric  is a  Finsler metric (on a manifold $M$)  of  the form
\begin{gather} \label{a0}  F(x,\xi ) =  \sqrt{g(\xi, \xi)}    +  \omega(\xi),
\end{gather}
where $g= g_{ij}$ is a Riemannian metric  and $\omega= \omega_i$ is an $1$-form.   The assumption that $F$ given by \eqref{a0} is indeed a Finsler metric is equivalent to the condition that the $g$-norm of $\omega$ is less than one.
Within the whole paper  we assume that all objects we consider are at least $C^2$-smooth.

     Two Finsler  metrics $F$ and $\bar F$
are {\it projectively equivalent}, if every forward-geodesic of~$F$   is a forward-geodesic of~$\bar F$.  Within our paper we  will always  assume that the dimension  of $M$ is at least two, since in dimension one all metrics are projectively equivalent.
 By  a  {\it  projective transformation}  of~$(M,F)$ we understand a  dif\/feomorphism  $a:M\to M$
such that pullback of~$F$ is projectively equivalent to~$F$.   The group of projective transformations  will be denoted by ${\rm Proj}(M,F)$; it is a Lie group of f\/inite dimension.

In \cite[Theorem~2]{iran} it was claimed that

\smallskip
{\it for  closed connected  $(n\ge 3)$-dimensional  Randers  manifold  $(M,F)$
of constant flag curvature $\dim({\rm Proj}(M,F))= n(n+2)$ or $\dim({\rm Proj}(M,F))\le  \frac{n(n+1)}{2}$ .}

\smallskip

The goal of the present  note  is to construct a counterexample to this statement: we construct a Riemannian manifold of constant sectional curvature   such that its group of projective transformations is $ (n-1)^2 + 1 = n^2 - 2 n + 2$ dimensional. Clearly, for   $n> 4 $
 we have
 \[
 \frac{n(n+1)}{2}< n^2 - 2 n + 2 = (n-1)^2 + 1< (n+1)^2 -1  =  n(n+2).
 \]
Of course, Riemannian manifolds of constant sectional curvature are Randers manifolds of constant f\/lag curvature.

\section{Example}

Let us f\/irst recall the group of the projective transformations of the standard sphere
\begin{gather} \label{1}
 S^n:=  \big\{(x_1,x_2, \dots,x_{n+1})\in \mathbb{ R}^{n+1}\mid  \ x_1^2+x_2^2+\cdots +x_{n+1}^2=1\big\}
\end{gather}
with the restriction of the standard Euclidean metric.
For every $A\in SL_{n+1}$, consider  the  dif\/feomorphism
\[
a:\ S^n\to S^n, \qquad
a:\ v\mapsto \frac{Av}{\|Av\|},
  \]
where  $Av$ is the product of the matrix $A$ and $v\in S^n $ viewed as a vector in $\mathbb{R}^{n+1}$  and $\|\cdot \|$ denotes the standard norm in $\mathbb{R}^{n+1}$.

The  dif\/feomorphism $a$ clearly  takes  geodesics to geodesics.
Indeed, the geodesics of  $g$ are great circles
(the intersections of planes that go through the origin with the
 sphere). Since  multiplication with $A$ is a linear mapping, it
 takes planes to  planes. Since the normalisation $w\mapsto
\frac{w}{\|w\|}$ takes punctured   planes to their intersections with the
sphere,  $a$ takes  great circles to great circles. Thus,  $a$ takes geodesic to geodesic.

It is wellknown (essentially, since the time of Beltrami \cite{Beltrami}, see also \cite{short}),  that all projective transformations of the sphere can be constructed by this procedure. Clearly, two matrices $A$ and $A'\in SL_{n+1}$ generate the same projective transformation if and only if $A'  = \pm  A$  (of course,
$A$ and $-A$ can  simultaneosly lie in $SL_{n+1}$ only if $n$ is odd).

Let us now construct  a quotient of the standard sphere modulo a (f\/inite,  freely acting) subgroup of the  isometry group, such that the   group of projective transformations of this quotient   has dimension $n^2 - 2 n + 2$.

Take $n\ge 2$ and consider the   orthogonal transformation of $\mathbb{R}^{n+1}$ whose matrix $B\in O_{n+1}$ is blockdiagonal such that  the f\/irst $2\times 2$-block equals   $\begin{pmatrix} 0 &  1
\\ -1 & 0\end{pmatrix}$ and the second $(n-1)\times (n-1)$-dimensional block equals  $-{\rm id} = {\rm diag}(-1,\dots,-1)$:
\begin{gather*}
B=  \begin{pmatrix} 0  & 1 &  &  & &  \\  -1 & 0 &  &  & & \\
 &  &-1   &  & & \\ &  &   &-1  & & \\ &  &   &  &\ddots  & \\ &  &   &  & & -1\end{pmatrix}.
 \end{gather*}

The restriction of this orthogonal transformation to the sphere will be denoted by $b$. Clearly, $b:S^n \to S^n$  is an involution  with no f\/ixed points preserving  the metric of the sphere.

Consider now the subgroup (isomorphic to $\mathbb{Z}_2$)  of the isometry group  of the sphere generated by the isometry $b$, and the  quotient $M:= S^n/  \mathbb{Z}_2$ with   the induced metric which we denote by~$g$. Let us  now show that the group of the projective transformations of $(M, g)$ has dimension $n^2 - 2 n + 2$.

We f\/irst construct a $(n^2 - 2 n + 2)$-dimensional  subgroup of ${\rm Proj}(M,g)$.   For every  $A\in SL_{n+1}$  such  the matrix of $A$ is blockdiagonal  such that the f\/irst $2\times 2$-block  has the form $\begin{pmatrix} \alpha  & \beta \\ -\beta & \alpha\end{pmatrix}, $ where $\alpha,\beta \in \mathbb{R} $ with $\alpha^2 + \beta^2 > 0$, and    the second $(n-1)\times (n-1)$-dimensional block is  an arbitrary matrix $\tilde A$ with determinant $\frac{1}{\alpha^2 + \beta^2}$,
\begin{gather} \label{A}
A =  \begin{pmatrix}
\alpha   & \beta  & 0   & \cdots   &  0 \\
-\beta  & \alpha  & 0  & \cdots   & 0   \\  0 &  0 & & &\\
\vdots  & \vdots   &  & \textrm{\huge $\tilde A$ } &  \\ 0 & 0  & & & \end{pmatrix}.   \end{gather}

  For every   such a matrix $A$, let us consider the projective transformations $a$
   of  $S^n$   given by \eqref{1}. Since the matrices $A$ and $B$ commute, the projective transformation $a$ commutes with the isometry $b$ and  induces a projective transformation of~$M$. Clearly, two matrices~$A$ and $A'$ of the form~\eqref{A} generate  the same  projective transformation of~$M$  if and only if
   $A'\in  \{ A, -A, BA, -BA\}$. We denote the group of such projective transformations
by ${\rm Proj}'$;  we evidently have
\[
 \dim({\rm Proj}')=  \dim(GL_{n-1})-1 + 2 = (n-1)^2- 1 + 2 =   n^2 - 2n + 2.
 \]

Let us show that the group ${\rm Proj}'$
  contains (and therefore, since it is connected,  coincides with) ${\rm Proj}_0(M, g)$,  where ${\rm Proj}_0(M, g)$ denotes
 the connected component of the group of  projective transformations of $(M,g)$
   containing the neutral element.

{\sloppy Indeed, every element ${\rm Proj}_0(M,g)$  from a  small neighborhood of the neutral element
   lies in a~certain one-parametric  subgroup of ${\rm Proj}_0(M,g)$. Now, one-parametric
    subgroups of ${\rm Proj}_0(M,g)$ are in one-to-one correspondence with projective vector f\/ields.  Since the pro\-per\-ty of a vector f\/ield to be projective is a local property, the lift of  every projective vector f\/ield  to~$S^n$ is   projective w.r.t.\ the standard metric on~$S^n$, so it generates a 1-parameter group of projective transformations on $S^n$ which we denote by $a_\tau$.
As we recalled above, for every projective transformation of $S^n$ there exists  a matrix $A$ such that  $a(x)= \frac{Ax}{\|Ax\|}$.
 We denote  by $\overset{\tau}{A}\in SL_{n+1}$  the  matrix  such that  $a_{\tau}(x)= \frac{\overset{\tau}{A}x}{\|\overset{\tau}{A} x\|}$.

 }

Since the   f\/low of such projective vector f\/ield   on $S^n$ evidently commutes with the action of $\mathbb{Z}_2$ generated  by $b$,
for any $\tau$   the transformation $a_\tau$ commutes  with $b$. Then, the matrices $\overset{\tau}{A}$ and $B$ also commute:
$B \overset{\tau}{A} - \overset{\tau}{A}  B =0$.
It is an easy exercise in a linear algebra to show that if  a  matrix $A\in SL_{n+1}$  commutes  with $B$, then it has
 the form \eqref{A}.   Thus, any projective transformation of $(M,g)$ from a small neighborhood of the neutral element of  ${\rm Proj}_0(M,g)$ lies in  ${\rm Proj}'$. Then, ${\rm Proj}_0(M,g)\subseteq  {\rm Proj}'$.
 Since  ${\rm Proj}'$ is evidently connected, ${\rm Proj}_0(M,g)= {\rm Proj}'$ as we claimed.

\section[The group of projective transformations of a closed Randers manifold]{The group of projective transformations\\ of a closed Randers manifold}

Above we constructed an example of a closed manifold $(M, g)$  whose dimension of the group of projective transformations contradicts the statement of \cite[Theorem~2]{iran}.   The universal cover of this manifold is the standard sphere with the standard metric, and has the dimension of the group of projective transformations equal to  $n(n+2)$, so the nature of our counterexample is in a certain sense topological. Of course, one can construct many other counterexamples of the same nature (i.e., such that the universal cover is the standard sphere but the group of the projective transformations has the dimension  less then the   dimension of the  group of the
projective transformations of  the standard sphere but is still suf\/f\/iciently big).

The following (recently proved) statement shows that such examples are essentially all possible counterexamples to \cite[Theorem 2]{iran}:

\smallskip

\noindent
{\bf Fact (Corollary~3 of \cite{projectiveranders})}.
 {\it Let $(M,F)$ be a closed connected Randers  manifold with~$F$ given by~\eqref{a0}.   Then, at least   one of the following possibilities holds:
  \begin{itemize}\itemsep=0pt
   \item[$(i)$] There exists a
   closed form $\hat \lambda$  such that    ${\rm Proj}_0(M,F)$ consists of isometries of
    the Finsler metric  $F(x,\xi)=  \sqrt{g(\xi,\xi) } + \omega(\xi) - \hat \lambda(\xi)  $, or
  \item[$(ii)$]  the form  $\omega$ is  closed  and  $g$ has constant positive  sectional curvature. \end{itemize}}

 In  the case  $(i)$ of Fact above,
i.e., when  there exists a closed form $\hat \lambda$  such that    ${\rm Proj}_0(M,F)$ consists of isometries of the     the Finsler metric  $F(x,\xi)=  \sqrt{g(\xi,\xi) } + \omega(\xi) - \hat \lambda(\xi) $, the group ${\rm Proj}(M,F)$   has dimension at most $\tfrac{n(n+1)}{2}$, see
\cite[Proposition 6.4]{troyanov}.  In the case $(ii)$ of Fact  above,
the group of projective transformations  of  $F$ coincides with the group of the projective transformations of~$g$. Now, by \cite{obata,archive}, the connected component ${\rm Proj}_0(M, g)$ of this group contains not only isometries,  only if~$g$  has constant positive sectional curvature. Thus, if the dimension of the group of the projective  transformations of a closed $(n\ge 2)$-dimensional
 Randers  manifold $(M, F)$  (with $F$ given by \eqref{a0})
  is greater than $\tfrac{n(n+1)}{2}$, then $\omega$ is a closed form and $g$ has constant positive sectional curvature. Then, the lift of~$F$   to the universal cover $\tilde M \cong S^n$
  of $M$  has the $n(n+2)$-dimensional   group of the projective transformations.

Note that, dif\/ferent from \cite{iran},  in the  Fact  above  we do not require that the metric has  constant f\/lag curvature~-- this is an additional condition on~$g$ and $\omega$  assumed   in   \cite[Theorem~2]{iran}. Let us explain that this assumption (together with the assumption that the manifold is closed) implies that the Randers Finsler metric such that the dimension of its projective group is at least
$\tfrac{n(n+1)}{2}$ is  actually a Riemannian metric (i.e., the form $\omega$ in the formula~\eqref{a0} is zero).

Indeed, if   the condition $(ii)$ of Fact above holds, the Finsler metric  $F$ is projectively f\/lat (i.e., in a neighborhood of every point one can f\/ind local coordinates such that the geodesics are straight lines).
If we additionally assume that the metric has constant f\/lag curvature, the  form~$\omega$ must be identically zero  by \cite[Section~7.3]{Bao} so the metric is Riemannian.

Now, the condition $(i)$ of Fact above  implies already that the metric~$F$ is a Riemannian metric in view of  \cite[Theorem~7.1]{troyanov} (for $n\ne 2,4$ this fact  also follows from  \cite{Deng,wang}).

\medskip
{\bf Acknowledgements.}
 We thank O.~Yakimova  for useful discussions and the Deutsche For\-schungsgemeinschaft (GK 1523) for  partial f\/inancial support.

\pdfbookmark[1]{References}{ref}
\LastPageEnding

\end{document}